\documentclass[11pt]{scrartcl}

\usepackage[T1]{fontenc}
\usepackage[utf8]{inputenc}
\usepackage{libertine}

\usepackage{amsmath,amssymb} 
\usepackage{ntheorem} 
\usepackage{graphicx} 
\usepackage{enumerate}
\usepackage{enumitem}
\usepackage{tikz}

\usepackage{microtype}
\usepackage{faktor}

\newcommand{\nn}{\mathbb N}
\newcommand{\zz}{\mathbb Z}
\newcommand{\qq}{\mathbb Q}
\newcommand{\rr}{\mathbb R}
\newcommand{\cc}{\mathbb C}

\newcommand\In{\subseteq}

\newcommand{\fakt}[2]{#1\,\raise1pt\hbox{\ensuremath{/}}#2} 
\newcommand{\mm}[1]{\;\!\!\!\!\mod{#1}}

\usepackage{caption}

\usetikzlibrary{calc}
\usetikzlibrary{intersections}

\newcommand{\oz}{\overline{z}}
\newcommand\oo[1]{\mathcal{O}_{\qq(#1)}}

\theoremheaderfont{\bfseries\sffamily}
\theorembodyfont{\rmfamily}
\theoremseparator{:}
\newtheorem*{bew}{Proof}
\newenvironment{beweis}{\begin{bew}}{\qed\end{bew}\vspace{2mm}}

\theoremheaderfont{\bfseries\sffamily}
\theorembodyfont{\rmfamily}
\theoremseparator{.}
\newcommand{\qed}{\hfill\ensuremath{\square}}

\newcounter{count}
\newtheorem{Satz}[count]{Theorem}
\newtheorem{Korollar}[count]{Corollary}
\newtheorem{Lemma}[count]{Lemma}

\newtheorem{Bemerkung}[count]{Remark}
\newenvironment{bem}{\begin{Bemerkung}}{\hfill$\sharp$\end{Bemerkung}}

\title{On origami rings}
\author{Dmitri Nedrenco\footnote{University of Wuerzburg, Department of Mathematics. dmitri.nedrenco@mathematik.uni-wuerzburg.de}}

\begin{document}
\maketitle

\begin{abstract}
 \noindent In  \cite{or} the authors investigate  the so called origami rings. Taking this paper as a starting point we find some further properties of origami rings.
\end{abstract}

 We work in the complex plane $\cc$ and identify it occasionally with $\rr^2$.  
Let $U$ be a set of ``directions" which are determined by complex numbers  $e^{i\alpha}$ for some angles $\alpha\in[0,2\pi)$. Two directions $e^{i\alpha}$ and $e^{i\beta}$ are equal iff $\alpha = \beta \!\mod\pi$.  Let $L_u(p)$ denote the line with the direction  $u$ through the point $p$, i.\;\!e. $L_u(p)=p+\rr u$.

Also, let $I_{u,v}(p,q)$ denote the intersection point of two lines $L_u(p)$ and $L_v(p)$ for two different directions $u,v$. We set $M_0:=\{0,1\}$ and define $M_j$ as the set of all intersection points $I_{u,v}(p,q)$ for $u\neq v$ and $p\neq q$ in which $u,v$ take on all values from $U$ and $p,q$ take on all values in $M_{j-1}$. Finally, we define $R(U):=\bigcup_{j\geq 0} M_j$. 

In \cite{or} the authors investigate the set $R(U)$ and they prove that $R(U)$ is a ring for every multiplicative semigroup $U$. We try to answer a part of a question asked in the cited work: is $R(U)$ a ring even if $U$ is not a semigroup?  

First we collect some properties of the points $I_{u,v}(p,q)$ which are all proved in \cite{or}.

\begin{Satz}
 Let $u=e^{i\alpha},v=e^{i\beta}\in U$ be two different directions and $p,q$ two different elements of $R(U)$. Moreover, let $s_{x,y}:= x\overline{y}-\overline{x}y$, where  $\bar{\cdot}$ means complex conjugation. The following statements then hold:
\begin{enumerate}\label{eigI}
 \item\label{eigI1}  $I_{u,v}(p,q) = \frac{s_{u,p}}{s_{u,v}}v + \frac{s_{v,q}}{s_{v,u}}u$.
\item\label{eigI2} $I_{u,v}(p,q)=I_{v,u}(q,p)$.
\item\label{eigI3} $I_{u,v}(p,q)= I_{u,v}(p,0)+I_{u,v}(0,q)$.
\item\label{eigI4} $I_{u,v}(p+q,0) = I_{u,v}(p,0)+I_{u,v}(q,0)$ und $I_{u,v}(rp,0)=rI_{u,v}(p,0)$ for all $r\in\rr$.
\item\label{eigI5} $I_{u,v}(0,1) = \frac{s_{v,1}}{s_{v,u}}u = \frac{\operatorname{Im}(v)}{\operatorname{Im}(v\overline{u})}u = \frac{\operatorname{sin}\beta}{\operatorname{sin}(\beta-\alpha)}e^{i\alpha}= \frac{1-v^2}{1-(\frac{v}{u})^2}$.
\end{enumerate}
Also $R(U)$ is an additive group.\hfill\ensuremath{\square}
\end{Satz}
In \cite{or} the authors pointed out that for some sets $U$ the set $R(U)$ is a ring even if $U$ is not a semigroup, for instance $R(0^\circ, 45^\circ,90^\circ)=\zz[i]$. There are also other obvious examples: \[R(0^\circ, 30^\circ,60^\circ)=\zz[e^{\frac{2\pi i}{3}}]=R(0^\circ, 60^\circ,120^\circ) \quad \text{and} \quad R(0^\circ, 45^\circ,135^\circ)=\zz[i].\]

\subsubsection*{Special case $U=\{1,u,v\}$}

One could conjecture that $R(U)$ is always a ring. However, we show that this is not the case by considering the ring structure of $R(U)$ for sets $U$ with three directions.

\begin{Satz}\label{ruz}
 Let $U=\{1,u,v\}$ with $u=e^{i\alpha}$ and $v=e^{i\beta}$ be given, where $0\neq \alpha\neq \beta \neq 0 \mod{\pi}$ holds. Moreover, let  
\[ z:=I_{u,v}(0,1) = \frac{s_{v,1}}{s_{v,u}}u.\]
Then we have $R(U) = \zz+z\zz$.
\end{Satz}

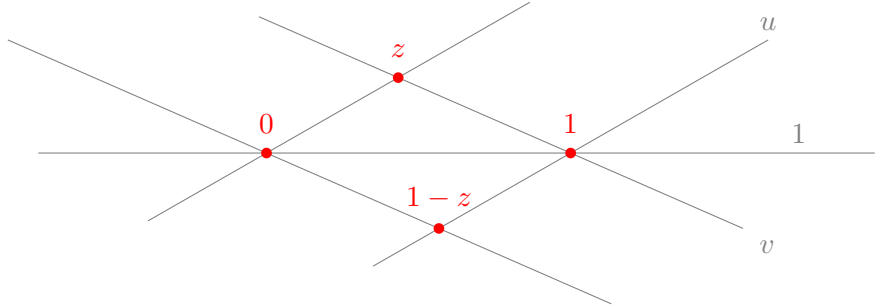
\begin{figure}[h!]
\begin{center}
 \begin{tikzpicture}
\begin{scope}[color=black!50]
\coordinate [label=above:$1$] () at (7,0);
\coordinate [label=below left:$v$] () at ([xshift=4cm]-20:3);
\coordinate [label=above:$u$] () at ([xshift=4cm]30:3);
\end{scope}

\draw[name path=xachse,black!50,thin] (-3,0)--(8,0);
\draw[name path=gerade1,black!50,thin] (210:1.8)--(0,0)--(30:4) [xshift=4cm] (210:3)--(0,0)--(30:3) ;
\draw[name path=gerade2,black!50,thin]($(30:2)!-2cm!(4,0)$)--($(4,0)!-1!(30:2)$)  [xshift=-4cm] ($(30:2)!-0.5!(4,0)$)--($(4,0)!-2!(30:2)$);
\fill[name intersections={of=gerade1 and gerade2},red] 
(intersection-1) circle (2pt) node[label=$z$] {} 
(intersection-2) circle (2pt) node[label=$1-z$] {};

\fill[name intersections={of=gerade1 and xachse},red] 
(intersection-1) circle (2pt) node[label=$0$] {};

\fill[name intersections={of=gerade2 and xachse},red] 
(intersection-1) circle (2pt) node[label=$1$] {};

 \end{tikzpicture}\end{center}
\caption{Construction of $M_1$.}
\end{figure}

\begin{beweis}

Since $1$ and $z$ belong to the additive group $R(U)$, we have $\zz+z\zz \In R(U)$. We prove the other inclusion by showing that $M_j\In \zz+z\zz$ via induction on $j$. This holds for $j=0$.
 Let  $s,t\in M_j$; then there exist $a,b\in\zz$ satisfying $s=a+bz$. Due to Theorem~\ref{eigI}\ref{eigI2} and \ref{eigI}\ref{eigI3} it suffices to show that $I_{x,y}(s,0)\in\zz+z\zz$ for $\{x,y\}\In U$.

We have to show that the following six points 
\[I_{u,v}(s,0),\; I_{v,u}(s,0),\; I_{u,1}(s,0),\; I_{v,1}(s,0),\; I_{1,u}(s,0),\; I_{1,v}(s,0)\]
belong to $\zz+z\zz.$ 

\noindent For this purpose we calculate
\begin{align*}
\tag{$\star$} z&=\tfrac{s_{v,1}}{s_{v,u}}u, \text{ here we have } s_{v,u}\neq 0 \text{ as } \alpha\neq \beta \!\!\!\mod{\pi},\; \tfrac{s_{v,1}}{s_{v,u}}\in\rr \text{ and}\\
\tag{$\star\star$} s_{u,z} &= u\overline{z}-\overline{u}z = u\overline{u}\big(\tfrac{s_{v,1}}{s_{v,u}}-\tfrac{s_{v,1}}{s_{v,u}}\big)=0.
\end{align*}
It is easy to \emph{see} how the points we are looking for are constructed. We provide an analytic proof by using Theorem \ref{eigI} although it is very helpful to draw a picture first in order to get better understanding of the calculations. 

\begin{align*}
\bullet\; I_{u,v}(s,0) = &~  I_{u,v}(a+bz,0) = aI_{u,v}(1,0)+bI_{u,v}(z,0) = a(1-z)\in\zz +z\zz,\\
& \text{ since } I_{u,v}(1,0) = 1-z \text{ and } I_{u,v}(z,0) = \frac{s_{u,z}}{s_{u,v}}v = 0 \text{ because of } (\star\star). \\[3mm]
\bullet\; I_{v,u}(s,0) = &~ aI_{v,u}(1,0)+bI_{v,u}(z,0) = (a+b)z\in\zz+z\zz,\\
&\text{ since } I_{v,u}(1,0)=z \text{ and } I_{v,u}(z,0) = I_{u,v}(0,z)=\frac{s_{v,z}}{s_{v,u}}u = z.\\[3mm]
\bullet\; I_{u,1}(s,0) = &~ aI_{u,1}(1,0)+bI_{u,1}(z,0) = a,\\
&\text{ since } I_{u,1}(1,0) = \frac{s_{u,1}}{s_{u,1}}\cdot 1 = 1 \text{ and } I_{u,1}(z,0) = \frac{s_{u,z}}{s_{u,1}}\cdot 1 = 0 \text{ because of } (\star\star).\\[3mm]
\bullet\; I_{v,1}(s,0) = &~  aI_{v,1}(1,0)+bI_{v,1}(z,0) = a+b,\\
&\text{ since } I_{v,1}(z,0)= \frac{s_{v,z}}{s_{v,1}}\cdot 1 = \frac{v\overline{z}-\overline{v}z}{s_{v,1}}\overset{(\star)}{=} \frac{(v\overline{u}-\overline{v}u)\cdot\ \frac{s_{v,1}}{s_{v,u}}}{s_{v,1}} = 1.\\[3mm]
\bullet\; I_{1,u}(s,0) = &~ aI_{1,u}(1,0)+bI_{1,u}(z,0) = bz
\text{ because of } I_{1,u}(1,0)= \frac{s_{1,1}}{s_{1,u}}u = 0\text{ and }\\& I_{1,u}(z,0)=\frac{s_{1,z}}{s_{1,u}}u = \frac{\overline{z}-z}{\overline{u}-u}u = \frac{\overline{u}-u}{\overline{u}-u}\cdot\frac{s_{v,1}}{s_{v,u}}u \overset{(\star)}{=} z.
\\[3mm]
\bullet\; I_{1,v}(s,0) = &~ aI_{1,v}(1,0)+bI_{1,v}(z,0) = b(z-1),
\text{ since } I_{1,v}(1,0) = 0 \text{ and }\\ 
 &I_{1,v}(z,0) = \frac{s_{1,z}}{s_{1,v}}v = \frac{\overline{z}-z}{\overline{v}-v}v = \frac{\overline{u}-u}{\overline{v}-v}\cdot \frac{s_{v,1}}{s_{v,u}}v  = \frac{\overline{u}v-uv}{\overline{v}-v}\cdot \frac{s_{v,1}}{s_{v,u}} = \frac{uv-\overline{u}v}{v\overline{u}-\overline{v}u}\\
& \text{but on the other hand }z-1 = \frac{s_{v,1}}{s_{v,u}}u - 1 = \frac{(v-\overline{v})u - v\overline{u}+\overline{v}u}{s_{v,u}}= \frac{uv-\overline{u}v}{v\overline{u}-\overline{v}u},\\
& \text{ so } I_{1,v}(z,0)=z-1.
\end{align*}
Thus, we dealt with all the cases, so the proof is complete.
\end{beweis}

\subsubsection*{General case for $|U|=3$}
In Theorem \ref{ruz} we assume that one of the directions in $U$ is determined by the angle of $0^\circ$. In fact this is no  loss of generality: 
Suppose $U=\{x,u,v\} $ is a set with three different directions and $M_0=\{0,1\}$, then we define $1'\!\!:=I_{x,v}(0,1)$. Because $GL_2\rr$ operates transitively on $\rr^2\backslash\{0\}$, we can transform  $M_0$ by a linear transformation to $\{0,1'\}$. 
Now we can in fact assume that one direction equals $1$. This means, up to a linear transformation, we proved that $R(U)$ is of the form  $\zz+z\zz$.

\subsubsection*{Ring structure of $R(1,u,v)$} Next, we want to clarify for which  directions $u=e^{i\alpha},v=e^{i\beta}$ different from $\pm 1$ where $\beta \neq \alpha \operatorname{mod}{\pi}$ the set $R(1,u,v)=\zz+z\zz$ is a ring. Since  \[ (a+bz)(c+dz) = ac + (ad+bc)z+bdz^2\]
we see that $\zz+z\zz$ is closed under multiplication iff $z^2$ lies in it. This amounts to the same as to say that the coefficients of the quadratic minimal polynomial of $z$ over $\rr$ are integers. This  is in turn the same as to say that $z$ is an algebraic number of degree $2$ (by the construction of $z$ it is not real). So for instance if  
$\alpha = \operatorname{arctan}\sqrt{7}$ and $\beta = \pi-\operatorname{arctan}\sqrt{7}$, then the set $R(1,e^{i\alpha}, e^{i\beta})$ is a ring, for it is $z^2=z-2$.

\begin{bem}\label{bem1i} Let $\alpha = \pm \frac{\pi}{2}$ i.\;\!e. $u=\pm i$. Since $u=i$ and $u=-i$ represent the same direction, we let without loss of generality $u=i$. The two directions $1$ and $i$ are perpendicular and\footnote{The following is not a restriction: we can consider $I_{i,v}(0,1)$ instead of $I_{v,i}(0,1)=1-z$.} let $z$ be $I_{i,v}(0,1)=ri$ for some $r\in\rr$. More precisely we have $r= \operatorname{tan}(\pi - \beta) = -\operatorname{tan}\beta\neq 0$ since $1,i,v$ are different directions. Hence the question whether $R(1,u,v)$ is a ring  is reduced to the question whether or not $z^2=-(\operatorname{tan}\beta)^2$ is an integer. From this we infer that  $R(U)$ is a ring for $U=\{1,i,v\}$ iff $\beta=\operatorname{arctan}\sqrt{d}$ where $d$ is a positive squarefree integer. In this case $v=exp({i\operatorname{arctan}\sqrt{d}})$ and $R(1,i,v)= \zz + i\tan\beta\,\zz = \zz+\sqrt{-d}\,\zz$.
\end{bem}

\begin{bem}\label{ruring}
The complex number $z=I_{u,v}(0,1)\in\cc, z\not\in\rr$ has the real minimal polynomial $(x-z)(x-\overline{z}) = x^2 -(z+\overline{z})x + z\overline{z}\in\rr[x]$. 
By the above considerations this means that $R(U)$ is a ring exactly when (cf. Theorem \ref{eigI}\;\!\ref{eigI5}\,):
\begin{align*}
k:&=z+\oz = 2\operatorname{Re}z =  2\cdot\frac{\operatorname{sin}\beta\operatorname{cos}\alpha}{\operatorname{sin}(\beta-\alpha)}\in\zz,\\
m:&=z\oz = |z|^2 = \frac{\operatorname{sin}^2\beta}{\operatorname{sin}^2(\beta-\alpha)}\in\zz.
\end{align*}
Thus, in particular $k = 2 \sqrt{m}\cos\alpha$ and $\cos^2\alpha = \frac{k^2}{4m}$ is a rational number.
 
Hence, if $R(U)$ is a ring, then necesssarily $\cos^2\alpha$ is a rational number. Since we could start with $u$ and $v$ interchanged it follows by symmetry that $\cos^2\beta\in \qq$. Obviously this property is not enough to guarantee that  $R(U)$ is a ring; this is clear by looking at the example  $R(0^\circ, 60^\circ,150^\circ)$, for $\cos 60^\circ = \frac{1}{2}$, but $k = 2\operatorname{Re}(z)=  \frac{1}{2} \not\in\zz$ (here  $z = \frac{1}{4}(1+i)$).

At least in the following sense the condition  $\cos\alpha\in\qq$ is sufficient (see figure~2): If $\cos\alpha = \frac{s}{t}$ and without loss of generality $s,t\in\nn$ and relatively prime, then we set $\operatorname{Re}z = s$ and $|z| =  t$. Therefore $k$ as well as $m$ are integers. With this choice we have $z=s+i\sqrt{t^2-s^2}$ and $\beta = \operatorname{arctan}(\frac{\sqrt{t^2-s^2}}{s-1})$. Hence, for every $\alpha$ with $\cos\alpha\in\qq$  we found  $\beta$ such that $R(1,u=e^{i\alpha},v=e^{i\beta})$ is a ring.
\end{bem}
\begin{Korollar}
 For infinitely many pairs $u,v$ of directions the set $R(1,u,v)$ is a ring.\hfill\qed
\end{Korollar}

 If $s$ and $t$ from the fraction $\frac{s}{t}$ in Remark \ref{ruring} are not relatively prime, we get other values of $\beta$ and $z$ (see figure 2), but since  $\cos\alpha = \frac{s}{t} =\frac{s\cdot \gamma}{t\cdot \gamma}$ where $\gamma\in\zz$ it holds that $z'=\gamma\cdot z$ so $\zz+z'\zz\In \zz+z\zz$ and we only obtain  a subring.

\begin{figure}[h!]\begin{center}
 \begin{tikzpicture}
\draw[name path=xachse,black!50,thin] (-2,0)--(7,0);
\begin{scope}[color=red]
\coordinate [label=above left:$z$] (z) at (5,3);
\coordinate [label=below:$0$] (0) at (0,0);
\coordinate [label=below:$1$] (1) at (2,0);
\coordinate [label=below right:${s=Re(z)}$] (s) at (5,0);
\end{scope}

\begin{scope}[color=black!50]
\coordinate [label=above:$\alpha$] (a) at (0.8,0);
\coordinate [label=above right:$\beta$] (b) at (2.6,0);
\coordinate [label=above:$1$] () at (-2,0);
\coordinate [label=below:$u$] () at ($(0)!-1.5cm!(z)$);
\coordinate [label=left:$v$] () at ($(1)!-1.6cm!(z)$); 
\end{scope}

\draw[thin,black!50] ($(0)!-1.5cm!(z)$)--($(z)!-1.8cm!(0)$)  node [sloped,midway,above, red]{$t=|z|$};
\draw[thin,black!50] ($(1)!-1.5cm!(z)$)--($(z)!-1.8cm!(1)$);
\draw[thick,black] (0)--(z)--(1) ;
\draw[thin, dashed,black!50] ($(z)!-1cm!(s)$)--($(s)!-1cm!(z)$);
\fill[red] (z) circle (2pt) (0) circle(2pt) (1) circle (2pt) (s) circle (2pt);

 \end{tikzpicture}\end{center}
\caption{Construction of $\beta$ for some given $\alpha$.}
\end{figure}
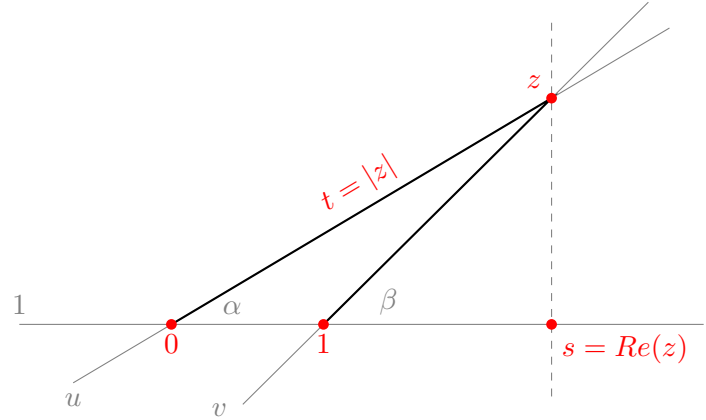

\subsubsection*{Quadratic number fields}
If we think of the minimal polynomial $x^2+kx+m$ of $z$ we see that $z=\frac{k}{2}+\frac{\sqrt{k^2-4m}}{2}$ (in what follows it does not matter whether we use $z$ or $\overline{z}$). Thus, the quadratic number field $\qq(z) = \faktor{\qq[x]}{(x^2+kx+m)\qq[x]}$ can be written as $\qq(z) = \qq(\sqrt{k^2-4m})$. Quadratic number fields are characterized as follows (cf. \cite[p.\,62, Proposition 3.1]{stewart}).

\begin{Satz}
 The quadratic number fields are exactly  the fields $\qq(\sqrt{d})$ for some squarefree integer $d$. \qed
\end{Satz}
Now we wish to know for which squarefree integers $d$ the field $\qq(\sqrt{d})$ equals $\qq(z)$ for a suitable choice  $k$ and $m$ (so actually of  $\alpha$ and $\beta$). 

\noindent Since $z\in\cc\backslash\rr$, such a $d$—if it exists—has to be negative, so we really look only at purely imaginary quadratic number fields. In Remark \ref{bem1i} we saw that if we choose $\alpha=\frac{\pi}{2}$ and $\beta=\operatorname{arctan}(\sqrt{d})$ where $d$ is a squarefree positive integer, then the equation $z^2=-(\operatorname{tan}\beta)^2 = -d$ holds. Hence for every squarefree positive integer $d$  there exist angles $\alpha$ and $\beta$ such that $\qq(z)=\qq(\sqrt{-d})$ is true.

\subsubsection*{Ring of algebraic numbers in $\qq(z)$}
Assume  that $R(1,u,v)=\zz+z\zz$ is a ring. Then  $z$ is an algebraic number of degree $2$. Since integers are for trivial reasons algebraic numbers, the ring $\mathcal{O}_{\qq(z)}$  of algebraic numbers  of the field  $\qq(z)$ contains the ring $\zz+z\zz$. 
We ask ourselves  when the equality holds. 

\begin{Lemma}\label{ganzalg}
 Let $d$ be a squarefree integer. Then the ring of algebraic numbers of $\qq(\sqrt{d})$ equals\footnote{For a proof see for instance \cite[p. 62, Theorem 3.1]{stewart}.}
\[ \oo{\sqrt{d}} =\begin{cases}
                  \zz+\sqrt{d}\,\zz& d\not\equiv 1 \!\!\mod 4\\
		  \zz + \frac{1+\sqrt{d}}{2}\,\zz& d\equiv 1\!\!\mod 4.
                 \end{cases}
\]
\end{Lemma}

\begin{Satz}\label{ganzalgring}
Let $1,u=e^{i\alpha},v=e^{i\beta}$ be pairwise different directions, let $z=I_{u,v}(0,1)$ and $k=z+\overline{z}$. Moreover let $R(U)=\zz+z\zz$ be a ring. In this case the ring $\oo{z}$ of algebraic integers of $\qq(z)$ equals $R(1,u,v)$ exactly in the following cases:
\begin{itemize}
 \item[] $k$ is odd and $(k\tan\alpha)^2$ is a squarefree positive integer; 
\item[{or:}] $k$ is even and $(\frac{k}{2}\tan\alpha)^2$ is a squarefree positive integer congruent to $2$ or $3$ modulo $4$;
\item[{or:}] $k=0$ and $\tan^2\!\beta$ is a squarefree positive integer congruent to $2$ or $3$ modulo $4$.
\end{itemize}

\end{Satz}

\begin{beweis}
We discuss the following cases: $k$ is odd, $0\neq k $ is even and  $k=0$.

In Remark \ref{bem1i} we dealt already with $k=0$: here $\oo{z}=R(U)=\zz+z\zz=\zz+\sqrt{-d}\,\zz$ iff $d$ is a squarefree positive integer congruent to $2$ or $3$ modulo $4$, therefore iff 
$\tan^2\!\beta$ is such a  $d$.

For the case $k\neq 0$ we calculate first the following:
\begin{equation}\tag{$\dagger$}
 \label{ktan}
k^2-4m=k^2(1-\frac{4m}{k^2})=k^2 (1-\frac{4m}{4m\cos^2\!\alpha})=k^2\cdot\frac{\cos^2\!\alpha-1}{\cos^2\!\alpha}=-(k\tan\alpha)^2.
\end{equation}
Let $k^2-4m = y^2 d $ where $d$ is a squarefree integer and $y\in\zz$.
Keeping Lemma \ref{ganzalg} in mind we ask when $\frac{1+\sqrt{d}}{2}\in\zz+z\zz$ is satisfied. Since $z = \frac{k}{2} + \frac{\sqrt{k^2-4m}}{2} = \frac{k}{2} + \frac{y\sqrt{d}}{2}$ this is the case exactly when  $y=\pm 1$ and $k$ is odd.

If $k$ is in fact  odd, then it follows that $k^2-4m \equiv 1\mod 4$ and so $y^2d \equiv k^2-4m \equiv 1\!\!\mod 4$ and  $y^2\equiv 1 \equiv d\!\mod 4$.  Therefore, for $k$ odd, the equation $\zz + z\zz = \oo{z}$ is true iff $(k\operatorname{tan}\alpha)^2$ is an odd squarefree integer, cf. (\ref{ktan}).

If $k$ is even and $\zz + z\zz = \oo{z}$ we infer that $y^2d \equiv k^2-4m \equiv 0\!\mod 4$ and as we have seen above $d\not\equiv1\mm 4$.
We check when $\sqrt{d}\in\zz + z\zz$. The fact that $d$ is squarefree enforces in $\sqrt{d} = a +b(\frac{k}{2}+\frac{y\sqrt{d}}{2})$ the conditions $y=\pm 2$ and $b = \pm 1$ (so  $a = \pm\frac{k}{2})$. Therefore,  $\sqrt{d}\in\zz + z\zz$ is true iff $k^2-4m = 4d$. That is why, for an even $k\neq 0$, the equation $\zz + z\zz = \oo{z}$ is satisfied  exactly when $(\frac{k}{2}\tan\alpha)^2$ is a squarefree positive integer congruent to $2$ or $3$ modulo $4$.
\end{beweis}
\noindent It does occur that $\zz+z\zz$ is a ring although $\zz+z\zz\neq \oo{z}$ holds: For $z=5+i\sqrt{56}$ we have $\cos\alpha = \frac{5}{9}$, thus  $\zz+z\zz$ is a ring according to Remark \ref{ruring}. However, $\tan\alpha = \frac{\sqrt{56}}{5}$ and $k=2\cdot 5$. Hence, we see that $\frac{k}{2}\tan\alpha = \sqrt{56}$ and $56$ is not squarefree. Thus, Theorem \ref{ganzalgring} says that $\zz+z\zz \varsubsetneq \oo{z}$.

\subsubsection*{The structure of $R(1,u,v,w)$} 
What can we say about $R(U)$ for a set $U$ consisting of $4$ different directions? As discussed above we can assume, up to a linear transformation, that the set $U$ equals $\{1,u,v,w\}$ for some different directions $1,u,v,w$. Moreover, we can choose $u,v,w$  such that $u=e^{i\alpha}, v= e^{i\beta}, w=e^{i\gamma}$ where $0<\alpha<\beta<\gamma<\pi$. Let \[p:=I_{u,w}(0,1) \quad \text{ and } \quad r:=I_{1,v}(p,0).\]
Then we have $r<p$ on the line $L_1(p)$ with its natural ordering.

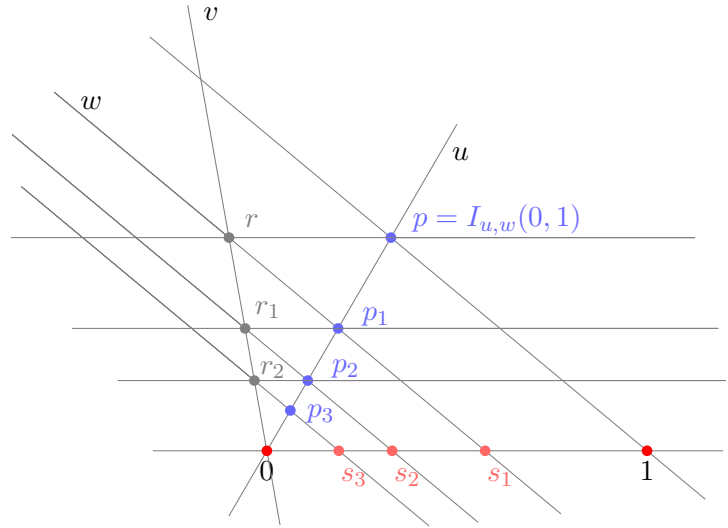
\begin{figure}[!b]
\vspace{-0mm}
\begin{center}
\begin{tikzpicture}
\coordinate [label=below:$0$] (0) at (0,0);
\coordinate [label=below:$1$] (1) at (5,0);
\coordinate  (t) at ([xshift=5cm]-40:1);

\draw[name path=xachse,black!50,thin] (-1.5,0)--(6,0);
\draw[name path=u,black!50,thin] (240:1)--(0,0)--(60:5);
\draw[name path=v,black!50,thin] (280:1)--(0,0)--(100:6);
\draw[name path=w1,black!50,thin] (t)--($(1)!-8.5cm!(t)$);

\fill[name intersections={of=u and w1, by={p}}][blue!60, every node/.style={right=7mm, below=0.5mm}]
(p) circle (2pt) node[label=${p=I_{u,w}(0,1)}$] {};
\draw[name path=x2, black!50,thin] (p)--++(4,0)--++(-9,0);

\fill[name intersections={of=x2 and v, by={r}}][black!50, every node/.style={right=1.5mm, below}]
(r) circle (2pt) node[label=$r$] {};
\draw[name path=w2, black!50,thin] (r)--++(140:3)--++(-40:8.7);

\fill[name intersections={of=w2 and u, by={p1}}][blue!60, every node/.style={right=2.5mm, below=0.5mm}]
(p1) circle (2pt) node[label=$p_1$] {};
\fill[name intersections={of=w2 and xachse, by={s1}}][red!60, every node/.style={right=1mm, below=3mm}]
(s1) circle (2pt) node[label=$s_1$] {};

\draw[name path=x3, black!50,thin] (p1)--++(5,0)--++(-8.5,0);
\fill[name intersections={of=x3 and v, by={r1}}][black!50, every node/.style={right=1.5mm, below}]
(r1) circle (2pt) node[label=$r_1$] {};

\draw[name path=w3, black!50,thin] (r1)--++(140:4)--++(-40:8);

\fill[name intersections={of=w3 and u, by={p2}}][blue!60, every node/.style={right=2.5mm, below=0.5mm}]
(p2) circle (2pt) node[label=$p_2$] {};
\fill[name intersections={of=w3 and xachse, by={s2}}][red!60, every node/.style={right=1mm, below=3mm}]
(s2) circle (2pt) node[label=$s_2$] {};

\draw[name path=x4, black!50,thin] (p2)--++(5.5,0)--++(-8,0);
\fill[name intersections={of=x4 and v, by={r2}}][black!50, every node/.style={right=1.2mm, below=0.5mm}]
(r2) circle (2pt) node[label=$r_2$] {};

\draw[name path=w4, black!50,thin] (r2)--++(140:4)--++(-40:7);

\fill[name intersections={of=w4 and u, by={p3}}][blue!60, every node/.style={right=2mm, below=1.5mm}]
(p3) circle (2pt) node[label=$p_3$] {};
\fill[name intersections={of=w4 and xachse, by={s3}}][red!60, every node/.style={right=1mm, below=3mm}]
(s3) circle (2pt) node[label=$s_3$] {};
\fill[red] (0) circle(2pt) (1) circle (2pt);

\coordinate [label={above:$w$}] () at ($(r)+(139:2.4)$);
\coordinate [label={right:$v$}] () at ($(r)+(99:3)$);
\coordinate [label={right:$u$}] () at ($(p)+(59:1.3)$);

\end{tikzpicture}
\end{center}
\caption{Construction of the sequences $(p_i)_i$ and $(s_i)_i$.}

\end{figure}

Since $L_w(r) \,||\, L_w(1)=L_w(p)$, it holds that $I_{1,w}(0,r)=L_1(0)\cap L_w(r) < 1$. But due to the choice $\gamma > \beta$, it follows that $I_{1,w}(p,0) < r$ on $L_1(p)$. Since $L_w(0)\, ||\, L_w(r)$ we have $I_{1,w}(0,r)=L_1(0)\cap L_w(r) > 0$. 

We infer that the triangle $0p1$ is similar to the triangle $0p_1s_1$, where \[ p_1:=I_{u,w}(0,r)\quad \text{ and } \quad s_1:=I_{1,w}(0,r) \quad\text{ with }  p_1\in(0,p), s_1\in (0,1) .\]
Using the point $p_1$ we construct another point on  $L_v(0)$, viz. $r_1:=I_{1,v}(p_1,0)$. Since $p_1\in(0,p)$ it follows for this point that $r_1\in(0,r)$.

\noindent With a similar argumentation as above we construct (see figure 3) \[p_2:=I_{u,w}(0,r_1) \quad \text{ and } \quad s_2:=I_{1,w}(0,r_1)\quad \text{ with } p_2\in(0,p_1), s_2\in (0,s_1).\]

\noindent Iteratively we construct the sequences $(p_i)_i$ and $(s_i)_i$ (as well as  the auxiliary sequence $(r_i)_i$) as follows:\vspace{-2pt}
\[
 p_i :=I_{u,w}(0,r_{i-1}),\quad s_i:=I_{1,w}(0,r_{i-1}),\quad r_{i-1}:=I_{1,v}(p_{i-1},0).
\]
Since $0$ is the only point on  $L_v(0)$ and $L_w(0)$, the points $p_i$ and $s_i$ are well defined. Furthermore the construction yields $p_i\in(0,p_{i-1})$ and $s_i\in(0,s_{i-1})$; therefore the triangles  $0p1$, $0p_is_i$ for $i=1,2,\ldots$ are all similar. In particular, due to compactness of the segments $(0,1)$ and $(0,p)$ each of the sequences $(p_i)_i$ and $(s_i)_i$ have a convergent subsequence.

\begin{Lemma}\label{abschluss}
 Let $X$  be a topologically closed subgroup of $\rr^2$, which is not contained in a line. Then there exists a basis $b,b'$ for $\rr^2$ such that 
\begin{itemize}
 \item[] $X = \rr b + \rr b' = \rr^2$ 
\item[{or}]  $X= \rr b + \zz b' $ 
\item[{or}] $X= \zz b + \zz b'$. \quad{For a proof see for instance} \cite[8.6, p. 83]{salzman}.\hfill\qed
\end{itemize}
\end{Lemma}

\begin{Korollar}
 If $U$ is a set containing at least $4$ different directions, then $R(U)$ is dense in $\cc$.
\end{Korollar}

\begin{beweis}
 \noindent Since $R(U)$ is an additive subgroup of $\rr^2$ (cf.\,Theorem \ref{eigI}), its closure has to be $\rr^2$, for the vectors $0p_1$ and $0s_1$ are linearly independent, cf. Lemma \ref{abschluss}. Therefore $R(U)$ is dense in $\rr^2$ resp. $\cc$.
\end{beweis}

\end{document}